\documentclass[11pt,twoside]{article}
\usepackage{mathbbold}

\usepackage{amsmath, amssymb, mathrsfs, graphicx}
\usepackage[all]{xy}

\def\az{\alpha}  
    \def\dz{\delta}
    \def\fz{\varphi}
  
\def\lz{\lambda} \def\mz{\mu}

\def\qd{\quad}
\def\qqd{\qquad}

\newcommand{\mathsym}[1]{{}}

\def\le{\leqslant}
\def\ge{\geqslant}

\font\cms=cmss9 scaled \magstep1

\def\nnd{\noindent}

\def\xmp{\nnd\bg{xmp1}}

\def\dexmp{\end{xmp1}}

\def\bg{\begin}
\def\be{\bg{equation}}
\def\de{\end{equation}}

\def\dear{\end{eqnarray}}
\def\lb{\label}
\def\ct{\cite}

\newcommand{\rf}[2]{[\ref{#1}; #2]}

\def\den{\end{enumerate}}

\begin{document}


\baselineskip 13pt

\thispagestyle{empty}
\renewcommand{\thefootnote}{\fnsymbol{footnote}}

\vspace*{-0.9in}
\noindent {Chapter ?? in \mbox{``Queueing Theory and Network Applications''},}\newline
\noindent {Springer£¬ 2017}

\begin{center}
{\bf\Large Trilogy on Computing Maximal Eigenpair}
\vskip.15in {Mu-Fa Chen}
\end{center}
\begin{center} (Beijing Normal University)\\
\vskip.1in June 8, 2017
\begin{center}

\end{center}
\end{center}
\vskip.1in

\markboth{\sc Mu-Fa Chen}{\sc Computing the maximal eigenpair}



\date{}



\nnd{\bf Abstract}\qd The eigenpair here means the twins consist of eigenvalue and its eigenvector.
This paper introduces the three steps of our study on computing the maximal eigenpair. In the first two steps, we construct efficient initials for a known but
dangerous algorithm, first for tridiagonal matrices and then for  irreducible matrices, having nonnegative off-diagonal elements. In the third step, we present two global algorithms
which are still efficient and work well for a quite large class of matrices, even complex for instance. \medskip

\nnd {\small 2000 {\it Mathematics Subject Classification}: 15A18, 65F15, 93E15}

\nnd {\small {\it Key words and phrases}. Maximal eigenpair, efficient initial, tridiagonal matrix, global algorithm.}

\bigskip

\section{Introduction}

This paper is a continuation of \ct{cmf17a}. For the reader's convenience,
we review shortly the first part of \ct{cmf17a}, especially the story of
the proportion of 1000 and 2 of iterations for two different algorithms.

The most famous result on the maximal eigenpair should be the Perron-Frobenius theorem. For nonnegative (pointwise) and irreducible $A$, if Trace\,$(A)$ $>0$, then the theorem says there exists uniquely a maximal eigenvalue $\rho(A)>0$ with positive left-eigenvector $u$ and positive right-eigenvector $g$ such that
      $$uA = \lz u, \qquad   A g = \lz g, \qquad   \lz= \rho (A).$$
These eigenvectors are also unique up to a constant.
Before going to the main body of the paper, let us make two remarks.

1) We need to study the right-eigenvector $g$ only. Otherwise, use the transpose $A^*$ instead of $A$.

2) The matrix\,$A$\,is required\,to\,be irreducible with nonnegative off-diagonal elements,
its diagonal elements can be arbitrary. Otherwise, use a shift $A + m I$ for large $m$:
\be (A+mI)g=\lz g\Longleftrightarrow Ag=(\lz-m)g,\de
their eigenvector remains the same but the maximal eigenvalues are shifted
to each other.

Consider the following matrix:
\be Q\!=\!\left(\!\begin{array}{ccccc}
-1^2 & 1^2 &0&0 &\cdots \\
1^2 & -1^2 - 2^2 & 2^2 &0 &\cdots \\
0& 2^2 & -2^2 - 3^2 & 3^2 &\cdots \\
\vdots &\vdots &\ddots &\ddots &\ddots \\
0& 0  & 0 &\; N^2 &-N^2-(N+1)^2
\end{array}\!\right)\!.\de
The main character of the matrix is the sequence $\{k^2\}$. The sum of each row equals zero except the last row. Actually, this matrix is truncated from the corresponding infinite one, in which case we have known that the maximal eigenvalue is $-1/4$ (refer to \rf{cmf10}{Example 3.6}).

\xmp\lb{t-01}\qd {\cms Let
$N=7$. Then the maximal eigenvalue is $-0.525268$ with eigenvector:
$$g\approx (55.878,\; 26.5271,\; 15.7059,\; 9.97983,\;
  6.43129,\; 4.0251,\; 2.2954,\; 1)^*,$$
where the vector $v^*=$ the transpose of $v$.
}\dexmp

We now want to practice the standard algorithms in matrix eigenvalue computation. The first method in computing the maximal eigenpair is the
{\it Power Iteration}, introduced in 1929.
Starting from a vector $v_0$ having a nonzero component in the direction of $g$, normalized with respect to a norm $\|\cdot\|$. At the $k$th step, iterate $v_k$ by the formula
\be v_{k}=\frac{A v_{k-1}}{\|A v_{k-1}\|},\quad {z_k}={\|Av_k\|},\qqd k\ge 1.\de
Then we have the convergence: $v_k\to g$ (first pointwise and then uniformly) and $z_k\to \rho(Q)$ as $k\to\infty$.
If we rewrite $v_k$ as
$$v_{k}=\frac{A^k v_{0}}{\|A^k v_{0}\|},$$
one sees where the name ``power'' comes from.
For our example, to use the Power Iteration, we adopt the $\ell^1$-norm and
choose $v_0={{\tilde v}_0}/{\|{\tilde v}_0\|}$, where
$${\tilde v}_0\!\!=\!\!(1,\,0.587624,\,0.426178,\,0.329975,\,0.260701,\,0.204394,0.153593,0.101142)^*\!\!.$$
This initial comes from a formula to be given in the next section. In Figure 1 below, the upper curve is $g$, the lower one is
modified from $\tilde v_0$, renormalized so that its last component becomes one. Clearly, these two functions are quite different, one may worry about the effectiveness of the choice of $v_0$. Anyhow, having the experience of computing its eigensystem, I expect to finish the computation in a few of seconds. Unexpectly, I got a difficult time to compute
{\begin{center}{\includegraphics[width=10.25cm,height=7.0cm]{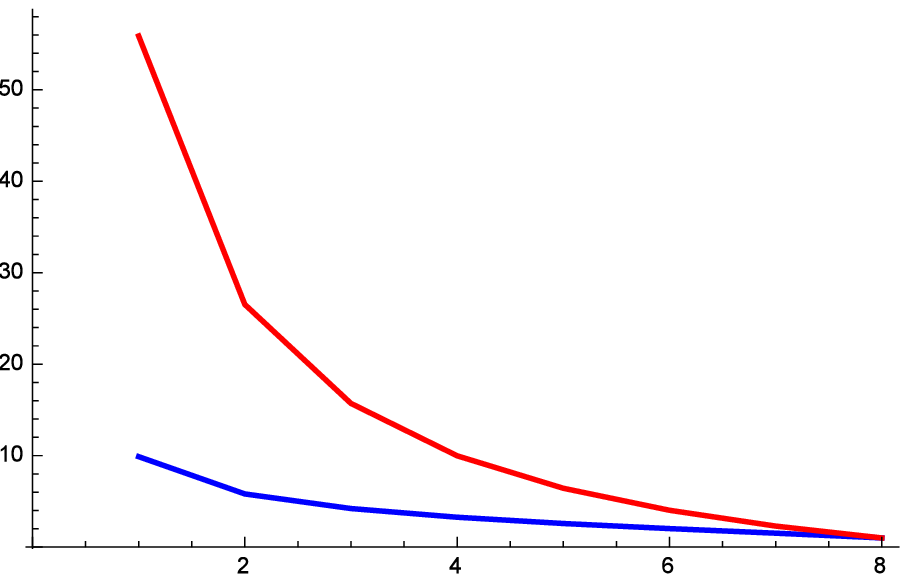}}\newline
{\vspace{-6truecm}\hspace{1truecm}{$$\begin{array}{ll}\text{The figure of $g$ and $v_0$}\\
\text{ }\end{array}$$}}
\end{center}}
\vspace{4truecm}
\begin{center} Figure 1: $g$ and $v_0$.\end{center}

\nnd the maximal eigenpair for this simple example. Altogether, I computed it for 180 times, not in one day, using 1000 iterations. The printed pdf-file of the outputs has 64 pages. Figure 2 gives us the outputs.

{\begin{center}{\includegraphics[width=12.25cm,height=9.0cm]{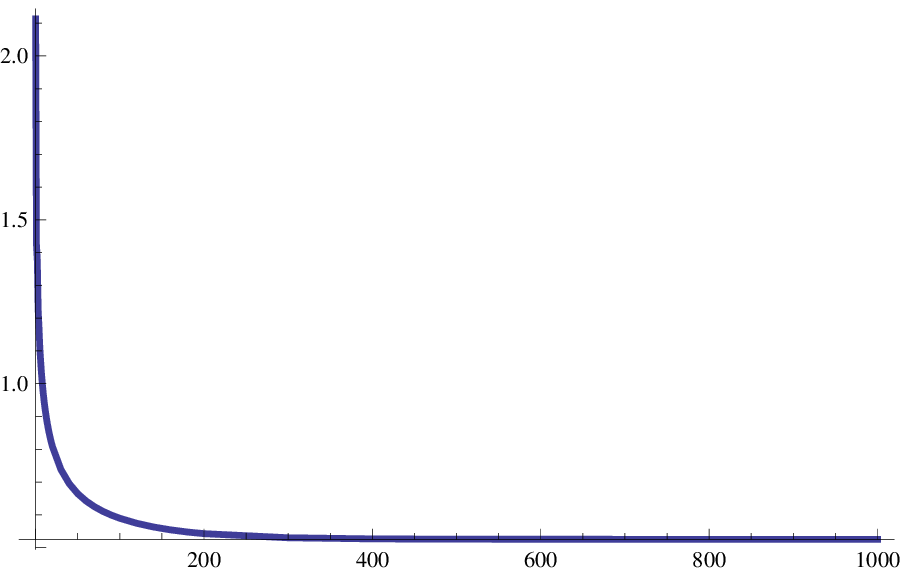}}\newline
{\vspace{-8truecm}\hspace{1truecm}{$$\begin{array}{ll}\text{The figure of } -z_k\\
\text{for }k=0, 1, \ldots, 1000.\end{array}$$}}
\end{center}}
\vspace{6truecm}
\begin{center} Figure 2: $-z_k$ for $k=0, 1, \ldots, 1000$.\end{center}

\nnd The figure shows that the convergence of $z_k$ goes quickly at the beginning of the iterations.
This means that our initial $v_0$ is good enough. Then the convergence goes very slow which means that the Power Iteration Algorithm converges very slowly.

Let us have a look at the convergence of the power iteration.
Suppose that the eigenvalues are all different for simplicity. Denote by
$(\lz_j, g_j)$ the eigenpairs with maximal one $(\lz_0, g_0)$. Write
$v_0=\sum_{j=0}^N c_j g_j$
for some constants $(c_j)$. Then $c_0\ne 0$ by assumption and
$$A^k v_0=\sum_{j=0}^N c_j\lz_j^k g_j
=c_0\lz_0^k \bigg[g_0+\sum_{j=1}^N \frac{c_j}{c_0}\bigg(\frac{\lz_j}{\lz_0}\bigg)^k g_j\bigg].
$$
Since $|\lz_j/\lz_0|<1$ for each $j\ge 1$ and $\|g_0\|=1$, we have
$$\frac{A^k v_0}{\|A^k v_0\|}= \frac{c_0}{|c_0|} g_0
+O\bigg(\bigg|\frac{\lz_1}{\lz_0}\bigg|^k\bigg)\qqd \text{as }k\to \infty,$$
where $|\lz_1|:=\max\{|\lz_j|: j>0\}.$
Since $|\lz_1/\lz_0|$ can be very closed to $1$,
this explains the reason why the convergence of the method can be very slow.

Before moving further, let us mention that the power method can be also used to compute the minimal eigenvalue $\lz_{\min}(A)$, simply replace $A$ by $A^{-1}$.  That is
the {\it Inverse Iteration} introduced in 1944:
\be v_k=\frac{A^{-1}v_{k-1}}{\|A^{-1}v_{k-1}\|}\Longleftrightarrow v_k=\frac{A^{-k}v_{0}}{\|A^{-k}v_{0}\|}.  \de
It is interesting to note that the equivalent assertion on the right-hand
side is exactly the the input-output method in economy.

To come back to compute the maximal $\rho(A)$ rather than $\lz_{\min}(A)$,
we add a shift $z$ to $A$: replacing $A$ by $A-z I$.
Actually, it is even better to replace the last one by $zI-A$ since we will
often use $z>\rho(A)$ rather than $z<\rho(A)$, the details will be explained
at the beginning of Section 4 below.
When $z$ is close enough to $\rho(A)$, the leading eigenvalue of $(z I-A)^{-1}$ becomes $(z-\rho(A))^{-1}$.
Furthermore, we can even use a variant shift $z_{k-1}I$ to accelerate the convergence speed.
Throughout this paper, we use varying shifts rather than a fixed one only.
Thus, we have arrived at
the second algorithm in computing the maximal eigenpair,
the {\it Rayleigh Quotient Iteration} (RQI), a variant of the {\it Inverse Iteration}. From now on, unless otherwise stated, we often use the $\ell^2$-norm.
Starting from an approximating pair $(z_0, v_0)$ of the maximal one $(\rho(A), g)$ with $v_0^*v_0=1$, use the following iteration.
 \be v_k=\frac{(z_{k-1}I-A)^{-1}v_{k-1}}{\|(z_{k-1}I-A)^{-1}v_{k-1}\|}, \qqd z_k=v_{k}^*Av_k,\qqd k\ge 1.\de
If $(z_0, v_0)$ is close enough to $(\rho(A), g)$, then
$$v_k\to g\qd \text{and}\qd z_k\to \rho(A) \qd \text{as }k\to\infty.$$
Since for each $k\ge 1$, $v_k^* v_k=1$, we have
$z_k={v_{k}^*Av_k}/{(v_{k}^* v_k)}.$
That is where the name ``Rayleigh Quotient'' comes from. Unless otherwise
stated, $z_0$ is setting to be $v_0^* A v_0$.

Having the hard time spent in the first algorithm, I wondered how many iterations are required using this algorithm.
Of course, I can no longer bear 1000 iterations. To be honest, I hope to finish the computation within 100 iterations. What happens now?

\xmp\qd{\cms For the same matrix $Q$ and ${\tilde v}_0$ as in Example \ref{t-01}, by RQI, we need two iterations only:}
$$z_1\approx -0.528215,\;\; z_2\approx -0.525268.$$
\dexmp

The result came to me, not enough to say surprisingly, I was shocked indeed.
This shows not only the power of the second method but also the effectiveness of my initial $v_0$. From the examples above, we have seen the story of
the proportion of 1000 and 2.

For simplicity, from now on, we often write $\lz_j:=\lz_j(-Q)$. In particular $\lz_0=-\rho(Q)>0$.
Instead of our previous $v_0$, we adopt the uniform distribution:
  $$v_0=(1, 1, 1, 1, 1, 1, 1, 1)^*/\sqrt{8}.$$
   This is somehow fair since we usually have no knowledge about $g$ in advance.
\xmp\lb{t-003}\qd{\cms Let $Q$ be the same as above. Use the uniform distribution $v_0$ and
set $z_0=v_0^*(-Q)v_0$. Then
$$\aligned
(z_1, z_2, z_3, {\pmb{z_4}}) &\approx  (4.78557,\;5.67061,\;5.91766,\;{\pmb{5.91867}}).\\
(\lz_0, \lz_1, {\pmb{\lz_2}}) &\approx  (0.525268,\,2.00758,\,{\pmb{5.91867}}).
\endaligned$$
}\dexmp
The computation becomes stable at the 4th iteration. Unfortunately, it is not what we want $\lz_0$ but $\lz_2$. In other words, the algorithm converges to a pitfall. Very often, there are $n-1$ pitfalls for a matrix having $n$ eigenvalues. This shows once again our initial ${\tilde v}_0$  is efficient
and the RQI is quite dangerous.

Hopefully, everyone here has heard the name {\it Google's PageRank}. In other words, the Google's search is based on the maximal left-eigenvector. On this topic, the book \ct{lm06} was published 11 years ago.
In this book, the Power Iteration is included but not the RQI. It should be clear that for PageRank, we need to consider not only large system, but also fast algorithm.

It may be the correct position to mention a part of the motivations for the present study.
\begin{itemize}  \setlength{\itemsep}{-0.8ex}
\item {Google's search}--PageRank.
\item {Input--output method} in economy. In this and the previous cases,  the computation of the maximal eigenvector is required.
\item {Stability speed} of stochastic systems. Here, for the
    stationary distribution of a Markov chain, we need to compute the eigenvector; and for the
    stability rate, we need to study the maximal (or the fist nontrivial) eigenvalue.
\item {Principal component analysis} for BigData. One choice is to study the so-called five-diagonal matrices. The second approach is using the maximal eigenvector to analysis the role played by the components, somehow similar to the PageRank.
\item For {image recognition}, one often uses Poisson or Toeplitz matrices, which are more or less the same as the Quasi-birth-death matrices studied in queueing theory. The discrete difference equations of elliptic partial differential equations are included in this class: the block-tridiagonal matrices.
\item The effectiveness of {random algorithm}, say Markov Chain Monte Carlo for instance, is described by the convergence speed. This is also related to the algorithms for machine learning.
\item As in the last item, a mathematical tool to describe the {phase transitions} is the first nontrivial eigenvalue (the next eigenpair in general). This is the original place where the author was attracted to the topic.
\end{itemize}
Since the wide range of the applications of the topic, there is a large number
of publications. The author is unable to present a carefully chosen list of references here, what instead are two random selected references: \ct{lm06} and \ct{sj15}.

Up to now, we have discussed only a small size $8\times 8 \,(N=7)$ matrix. How about large $N$? In computational mathematics, one often expects the number
of iterations grows in a polynomial way $N^{\az}$ for $\az$ greater or equal to 1. In our efficient  case, since $2=8^{1/3}$, we expect to have $10000^{1/3}\approx 22$ iterations
for $N\!+\!1\!=\! 10^4$. The next table subverts completely my imagination.

\begin{center}{\rm {\bf Table 1} \quad Comparison of RQI for different $N$}\end{center}
\vspace{-0.5truecm}

\begin{center}
{\begin{tabular}{|c|c|c|c|c|}
   \hline
$\pmb{N+1}$ & $\pmb{z_0}$ & $\pmb{z_1}$ & {$\pmb{z_2=\lz_0}$} & {\bf upper/lower} \\
  \hline\hline
$8$ & 0.523309 & {0.525268} & 0.525268 & $1\!+\!10^{-11}$\\
  \hline
$100$ & 0.387333 & 0.376393 & 0.376383 & $1\!+\!10^{-8}$\\
   \hline
$500$ & 0.349147 & 0.338342 & 0.338329 & $1\!+\!10^{-7}$\\
   \hline
\!$1000$\! & 0.338027 & 0.327254 & 0.32724 & $1\!+\!10^{-7}$\\
   \hline
\!$5000$\! & 0.319895 & 0.30855 & 0.308529 & $1\!+\!10^{-7}$\\
   \hline
\!$7500$\! & 0.316529 & 0.304942 & 0.304918 & $1\!+\!10^{-7}$\\
   \hline
$10^4$ & 0.31437 & 0.302586 & 0.302561 & $1\!+\!10^{-7}$\\
   \hline
\end{tabular}
}\end{center}
Here $z_0$ is defined by
$$z_0={7}/({8\dz_1})+v_0^*(-Q)v_0/8,$$
where ${v}_0$ and $\dz_1$ are computed by our general formulas to be defined in the next section.
We compute the matrices of order $8, 100,\ldots , 10^4$ by using MatLab in a notebook, in no more than 30 seconds, the iterations finish at the second step. This means that the outputs starting from $z_2$ are the same and coincide with $\lz_0$. See the first row for instance, which becomes stable at the first step indeed. We do not believe such a result for some days, so we checked it in different ways. First, since $\lz_0=1/4$ when $N=\infty$, the answers of $\lz_0$ given
in the fourth column are reasonable. More essentially, by using the output $v_2$, we can deduce upper and lower bounds of $\lz_0$ (using \rf{cmf10}{Theorem 2.4\,(3)}), and then the ratio upper/ lower is presented in the last column. In each case, the algorithm is significant up to 6 digits. For the large scale matrices here and in \ref{s-03}, the computations are completed by Yue-Shuang Li.

\section{Efficient initials: tridiagonal case}

   It is the position to write down the formulas of $v_0$ and $\dz_1$. Then our initial $z_0$ used in Table 1 is a little modification of $\dz_1^{-1}$: a convex combination of $\dz_1^{-1}$ and $v_0^*(-Q)v_0$.

   Let us consider the tridiagonal matrix (cf. \rf{cmf16}{\S 3} and \rf{cmf17c}{\S 4.4}). Fix $N\ge 1$, denote by
       $E=\{0, 1, \ldots, N\}$
the set of indices. By a shift if necessary, we may reduce $A$ to $Q$ with negative diagonals:
$Q^c=A-m I$, $ m:=\max_{i\in E} \sum_{j\in E} a_{ij},$
$$Q^c\!=\!\left(\!\begin{array}{ccccc}
-b_0-c_0 & b_0 &0&0 &\cdots \\
a_1 & -a_1 - b_1-c_1 & b_1 &0 &\cdots \\
0& a_2 & -a_2 - b_2-c_2 & b_2 &\cdots \\
\vdots &\vdots &\ddots &\ddots &\ddots \\
0& 0  & 0 &\; a_N &-a_N-c_N
\end{array}\!\right)\!.$$
Thus, we have three sequences $\{a_i>0\}$, $\{b_i>0\}$, and $\{c_i\ge 0\}$.
Our main assumption here is that the first two sequences are positive and $c_i\not\equiv 0$. In order to define our initials, we need three new sequences, $\{h_k\}$,
$ \{\mu_k\}$, and $\{\fz_k\}$.

First, we define the sequence $\{h_k\}$:
\be h_0=1,\;\; h_n=h_{n-1}r_{n-1},\qqd 1\le n \le N;\de
here we need another sequence $\{r_k\}$:
$$r_0=1+\frac{c_0}{b_0},\;\;
r_n=1+\frac{a_n+c_n}{b_n}-\frac{a_n}{b_n r_{n-1}},\qqd 1\le n <N.$$
Here and in what follows, our iterations are often of one-step.
Note that if $c_k= 0$ for every $k<N$, then we do not need the sequence $\{h_k\}$, simply set $h_k\equiv 1$.
An easier way to remember this $(h_i)$ is as follows. It is nearly harmonic of $Q^c$ except at the last point $N$:
\be Q^{c\,\setminus \text{\rm the last row}}h=0,\de
where $B^{\setminus \text{\rm the last row}}$ means the matrix modified from $B$ by removing its last low.

We now use $H$-transform, it is designed to remove the sequence $(c_i)$:
$${\widetilde Q}=\text{\rm Diag}(h_i)^{-1}Q^c\,\text{\rm Diag}(h_i).$$
Then
$$\aligned
{\widetilde Q}&=\left(\!\begin{array}{ccccc}
-b_0 & b_0 &0&0 &\cdots \\
a_1 &  -a_1 - b_1 & b_1 &0 &\cdots \\
0& a_2 & -a_2 - b_2 & b_2 &\cdots \\
\vdots &\vdots &\ddots &\ddots &\ddots \\
0& 0& \qquad 0 &\quad a_N^{}\;\; &-a_N^{}\!-\!c_N^{}
\end{array}\!\!\right)
\endaligned$$
for some modified $\{a_i>0\}$, $\{b_i>0\}$, and $c_N^{}> 0$.
Of course, $Q^c$ and ${\widetilde Q}$ have the same spectrum. In particular,
under the $H$-transform,
$$(\lz_{\min}(-Q^c),\; g)\to \big(\lz_{\min}\big(-{\widetilde Q}\big)=\lz_{\min}(-Q^c),\; \text{\rm Diag}(h_i)^{-1}g\big).$$
From now on, for simplicity, we denote by $Q$ the matrix replacing $c_N^{}$
by $b_N^{}$ in ${\widetilde Q}$.

Next, we define the second sequence $\{\mu_k\}$:
\be \mu_0=1,\;\;\mu_n=\mu_{n-1}\frac{b_{n-1}}{a_n},\qqd 1\le n\le  N.\de
And then define the third one $\{\fz_k\}$ as follows:
\be \fz_n=\sum_{k= n}^N \frac{1}{\mu_k b_k}, \qqd 0\le n\le N.\de

We are now ready to define $v_0$ and $\dz_1$ (or $z_0$) using the sequences
$(\mu_i)$ and $(\fz_i)$.
\begin{align}
&{{\tilde v}_0(i)\!=\!\sqrt{\fz_i}},\;i\le N; \qqd {v_0}\!=\!{\tilde v}_0/\|{\tilde v}_0\|;\qd \|\cdot\|:=\|\cdot\|_{L^2(\mu)}\\
&{\dz_1}\!=\!\max_{0\le n\le N} \bigg[\sqrt{\fz_n} \sum_{k=0}^n \mu_k \sqrt{\fz_k}
  +\!\frac{1}{\sqrt{\fz_n}}\sum_{n+1\le j \le N}\!\!\mu_j\fz_j^{3/2}\bigg]\!\!=:\!{z_0^{-1}}\!\!
 \end{align}
with a convention $\sum_{\emptyset}=0$.

   Finally, having constructed the initials $(v_0, z_0)$, the RQI goes as follows.
Solve $w_k$:
\be (-Q-z_{k-1}I) w_k=v_{k-1}, \qqd k\ge 1; \lb{11}\de
and define
$$v_k={w_{k}}/\!{\|w_k}\|,\qqd z_k=(v_{k},\, -Q\, v_k)_{L^2(\mu)}.$$
Then
$$v_k\to g \qd \text{and}\qd z_k\to \lz_0\qqd \text{as } k\to\infty.$$

Before moving further, let us mention that there is an explicit representation of
the solution $(w_i)$ to equation (\ref{11}). Assume that we are given $v:=v_{k-1}$ and
$z:=z_{k-1}$. Set
\be M_{sj}=\mu_j\sum_{k=j}^s \frac {1}{\mu_k b_k},\qqd 0\le j\le s\le N.  \lb{10}\de
Define two independent sequences $\{A(s)\}$ and $\{B(s)\}$, recurrently:
\be{\begin{cases}
A(s)= -\sum_{0\le j\le s-1} M_{s-1, j}\big(v(j) + z A(j)\big),\\
B(s)=1-z \sum_{0\le j\le s-1} M_{s-1, j} B(j),\qqd\qd 0\le s\le N.
\end{cases}} \de
Set
\be x=\frac{\sum_{j=0}^N \mu_j\big(v(j)+zA(j)\big)-\mz_N^{} b_N A(N)}{\mz_N^{} b_N B(N)-z \sum_{j=0}^N \mu_j B(j)}. \lb{12}  \de
Then the required solution $w_k:=\{w(s): s\in E\}$ can be expressed as $w(s)=A(s)+ x B(s)\,(s\in E)$.

To finish the algorithm, we return to the estimates of $\big(\lz_{\min}(-Q^c), g(Q^c)\big)$ ($g(Q^c)=g(-Q^c)$) or further $(\rho(A), g(A))$ if necessary, where $g(A)$, for instance, denotes the maximal eigenvector of $A$. Suppose that the iterations are stopped at $k=k_0$ and set $(\bar z, \bar v)=\big(z_{k_0}^{}, v_{k_0}^{}\big)$ for simplicity. Then, we have
$$\big(\lz_{\min}\big(-{Q^c}\big),\; \text{\rm Diag}(h_i)^{-1}g(Q^c)\big)=
\big(\lz_{\min}\big(-{\widetilde Q}\big),\; g\big({\widetilde Q}\big) \big)
\approx (\bar z, \bar v),$$
and so
\be \big(\lz_{\min} (-Q^c),\; g(Q^c)\big)
\approx \big(\bar z,\; \text{\rm Diag}(h_i)\,{\bar v}\big).\de
Because $\lz_{\min} (-Q^c)=m-\rho(A)$, we obtain
\be (\rho(A),\; g(A))\approx \big(m-{\bar z},\; \text{\rm Diag}(h_i)\,{\bar v}\big).\de

Now, the question is the possibility from the tridiagonal case to the general one.

\section{Efficient initials: the general case (\rf{cmf16}{\S 4.2} and \rf{cmf17c}{\S 4.5})}

When we first look at the question just mentioned, it seems quite a long distance to go from the special tridiagonal case to the general one. However, in the eigenvalue
computation theory, there is the so-called Lanczos tridiagonalization procedure
to handle the job, as discussed in \rf{cmf16}{Appendix of \S 3}. Nevertheless, what we adopted in \rf{cmf16}{\S 4} is a completely different approach. Here is our
main idea. Note that the initials $v_0$ and $\dz_1$ constructed in the last section are explicitly expressed by the new sequences. In other words,
we have used three new sequences $\{h_k\}$, $ \{\mu_k\}$, and $\{\fz_k\}$ instead of the original three $\{a_i\}$, $\{b_i\}$, and $\{c_i\}$ to describe our initials.
Very fortunately, the former three sequences do have clearly the
probabilistic meaning, which then leads us a way to go to the general setup.
Shortly, we construct these sequences by solving three linear equations
(usually, we do not have explicit solution in such a general setup). Then
use them to construct the initials and further apply the RQI-algorithm.

Let $A=(a_{ij}: i,j\in E)$ be the same as given at the beginning
of the paper. Set $A_i=\sum_{j\in E} a_{ij}$ and define
$$Q^c=A-\Big(\max_{i\in E} A_i\Big) I.
$$
We can now state the probabilistic/analytic meaning of the required three
sequences $(h_i)$, $(\mu_i)$, and $(\fz_i)$.
\bg{itemize}  \setlength{\itemsep}{-0.8ex}
\item $(h_i)$ is the harmonic function of $Q^c$ except at the right endpoint $N$, as mentioned in the last section.
\item $(\mu_i)$ is the invariant measure (stationary distribution) of the matrix
      $Q^c$ removing the sequence $(c_i)$.
\item $(\fz_i)$ is the tail related to the transiency series, refer to \rf{cmf16}{Lemma 24 and its proof}.
\end{itemize}

We now begin with our construction. Let
${h}=(h_0, h_1,\ldots, h_N)^*$ (with $h_0=1$) solve the equation
$$Q^{c\;\setminus \text{\rm the last row}}h = 0$$
and define
$${\widetilde Q}=\text{Diag}(h_i)^{-1}Q^c\,\text{Diag}(h_i).$$
Then for which we have
$$c_0^{}=\ldots=c_{N-1}^{}=0,\qd c_N^{}=:q_{N, N+1}^{}>0.$$
This is very much similar to the tridiagonal case.

Next, set $Q={\widetilde Q}$. Let ${\fz}=(\fz_0, \fz_1,\ldots, \fz_N)^*$ (with $\fz_0=1)$ solve the equation
$${\fz^{\setminus \text{\rm the first row}}=P^{\setminus \text{\rm the first row}}\,\fz},$$
where
$$P=\text{\rm Diag}\big((-q_{ii})^{-1}\big)Q+I.
$$

Thirdly, assume that ${\mu}:=(\mu_0, \mu_1, \ldots, \mu_N)$ with $\mu_0=1$ solves the equation
$${{Q}^{*\,\setminus \text{\rm the last row}}\mu^*=0}.$$

Having these sequences at hand, we can define the initials
$${\tilde v}_0(i)= \sqrt{\fz_i}, \;\;i\!\le\! N;\qd
v_0={\tilde v}_0/\|{\tilde v}_0\|_{\mu};\qd z_0=(v_{0}, -Q v_0)_{\mu}.
$$
Then, go to the RQI as usual. For $k\ge 1$, let $w_k$ solve the equation
$$(-Q-z_{k-1} I)w_{k}=v_{k-1}$$
and set
$$v_k={w_{k}}/\|w_k\|_{\mu},\qqd
  z_k=(v_{k}, -Q v_k)_{\mu}.$$
Then we often have $(z_k, v_k)\to (\lz_0, g)$ as $k\to\infty$.

We remark that there is an alternative choice (more safe) of $z_0$:
$$z_0^{-1}=\frac{1}{1-\fz_1}\max_{0\le n\le N}\bigg[\sqrt{\fz_n}
\sum_{k=0}^n \mu_k \sqrt{\fz_k}
\!+\!\frac{1}{\sqrt{\fz_n}} \sum_{n+1\le j\le N}\mu_j \fz_j^{3/2}\bigg]$$
which is almost a copy of the one used in the last section.

The procedure for returning to the estimates of $\big(\lz_{\min}(-Q^c), g(Q^c)\big)$ or further $(\rho(A), g(A))$ is very much the same as in the
last section.

To conclude this section, we introduce two examples to illustrate the
efficiency of the extended initials for tridiagonally dominant matrices.
The next two examples were computed by Xu Zhu, a master student in Shanghai.

\xmp[Block-tridiagonal matrix]\lb{t-02}\qd{\cms
Consider the matrix
$$Q=\left(\!\!\begin{array}{ccccc}
\!A_0 & B_0 &0&0 &\cdots \\
C_1 &A_1 & B_1 &0 &\cdots \\
0& C_2 & A_2 & B_2 &\cdots \\
\vdots &\vdots &\ddots &\ddots &\ddots \\
0& 0  & 0 &\quad C_N & A_N
\end{array}\right),$$
where $A_k, B_k, C_k$ are $40\times 40$-matrices,
$B$'s and $C$'s are identity matrices, and $A$'s are tridiagonal
matrices. For this model, two iterations are enough to arrive at
the required results (Table 2).
\begin{center}{\rm {\bf Table 2} \quad Outputs for Poisson matrix}\end{center}
\vspace{-0.6truecm}
\bg{center} {\begin{tabular}{|c|c|c|c|}
   \hline
$\pmb{\!N\!+\!1\!}$ & $\pmb{z_0}$ & $\pmb{z_1}$ &\pmb {$z_2=\lz_0$} \\
  \hline\hline
$1600$ & $7.985026$ & $7.988219$ & $7.988263$ \\
   \hline
 $3600$ & $7.993232$ & $7.994676$ & $7.994696$ \\
   \hline
 {$6400$} & $7.996161$ & $7.988256$ & $7.987972$ \\
   \hline
\end{tabular}}\end{center}
}\dexmp

\xmp[Toeplitz matrix]\lb{t-02}{\cms
Consider the matrix
$$A\!=\!\left(\!\!\begin{array}{cccccc}
1 & {2} &3&\cdots& n\!-\!1 & n \\
2 & 1 &{2}& \cdots &n\!-\!2 &n\!-\!1\\
\vdots &\vdots&\vdots &\ddots &\vdots &\vdots \\
{n\!-\!1}& n\!-\!2 & n\!-\!3&\cdots & 1 & {2} \\
n& {n\!-\!1} &n\!-\!2 & \cdots & 2 & 1
\end{array}\right).$$
 For this model, three iterations are enough to arrive at
the required results (Table 3).
\begin{center}{\rm {\bf Table 3} \quad Outputs for Toeplitz matrix}\end{center}
\vspace{-0.6truecm}
\begin{center}{\begin{tabular}{|c|c|c|c|c|}
 \hline
$\pmb{\!N\!+\!1\!}$ & $\pmb{z_0{\times 10^6}}$ & $\pmb{z_1{\times 10^6}}$ & $\pmb{z_2{\times 10^6}}$&\pmb{$z_3\!=\!\lz_0$}\\
  \hline\hline
$1600$ & 0.156992 & 0.451326 & 0.390252 &\!\!$0.389890$\! \\
   \hline
 $3600$ & 0.157398 & 2.30731 & 1.97816 &\!\!$1.97591$\! \\
   \hline
 {$6400$} & 0.157450 & 7.32791 & 6.25506 &\!\!$6.24718$\\
   \hline
\end{tabular}}\end{center}
}\dexmp

As mentioned before, the extended algorithm should be powerful for the tridiagonally dominant matrices. How about more general case?
Two questions are often asked to me by specialists in computational mathematics: do you allow more negative off-diagonal elements? How about complex matrices?
My answer is: they are too far away from me, since those matrices can not
be a generator of a Markov chain, I do not have a tool to handle them.
Alternatively, I have studied some more general matrices than the tridiagonal ones:
the block-tridiagonal matrices, the lower triangular plus upper-diagonal, the upper triangular plus lower-diagonal, and so on.
Certainly, we can do a lot case by case, but this seems still a long way
to achieve a global algorithm. So we do need a different idea.

\section{Global algorithms}\lb{s-03}

Several months ago, AlphaGo came to my attention. From which I learnt the subject of machine learning. After some days, I suddenly thought, since we are doing the computational mathematics, why can not let the computer help us to find a high efficiency initial value? Why can not we leave this hard task to the computer? If so, then we can start from a relatively simple and common initial value, let the computer help us to gradually improve it.

The first step is easy, simply choose the uniform distribution
as our initial $v_0$:
$$v_0=(1, 1, \cdots, 1)^*/\sqrt{N+1}.$$
As mentioned before, this initial vector is fair and universal. One may
feel strange at the first look at ``global'' in the title of this section.
However, with this universal $v_0$, the power iteration is already a global algorithm.
Unfortunately, the convergence of this method is too slow, and hence is
often not practical. To quicken the speed, we should add a shift which
now has a very heavy duty for our algorithm. The main trouble is that
the usual Rayleigh quotient
${v_0^*Av_0}/({v_0^*v_0})$
can not be used as $z_0$, otherwise, it will often lead to a pitfall, as illustrated by Example \ref{t-003}. The
main reason is that our $v_0$ is too rough and so $z_0$ deduced from it
is also too rough. Now, how to choose $z_0$ and further $z_n$?

Clearly, for avoiding the pitfalls, we have to choose $z_0$ from the outside
of the spectrum of $A$ (denoted by Sp$(A)$), and as close to $\rho(A)$ as possible to quicken the convergence speed. For nonnegative $A$, Sp$(A)$ is
located in a circle with radius $\rho(A)$ in the complex plane. Thus, the safe
region should be on the outside of Sp$(A)$. Since $\rho(A)$ is located at the boundary on the right-hand side of the circle, the effective area should be on the real axis on the right-hand side of, but a little away from, $\rho(A)$.

{\begin{center}{\includegraphics[width=6.25cm,height=3.5cm]{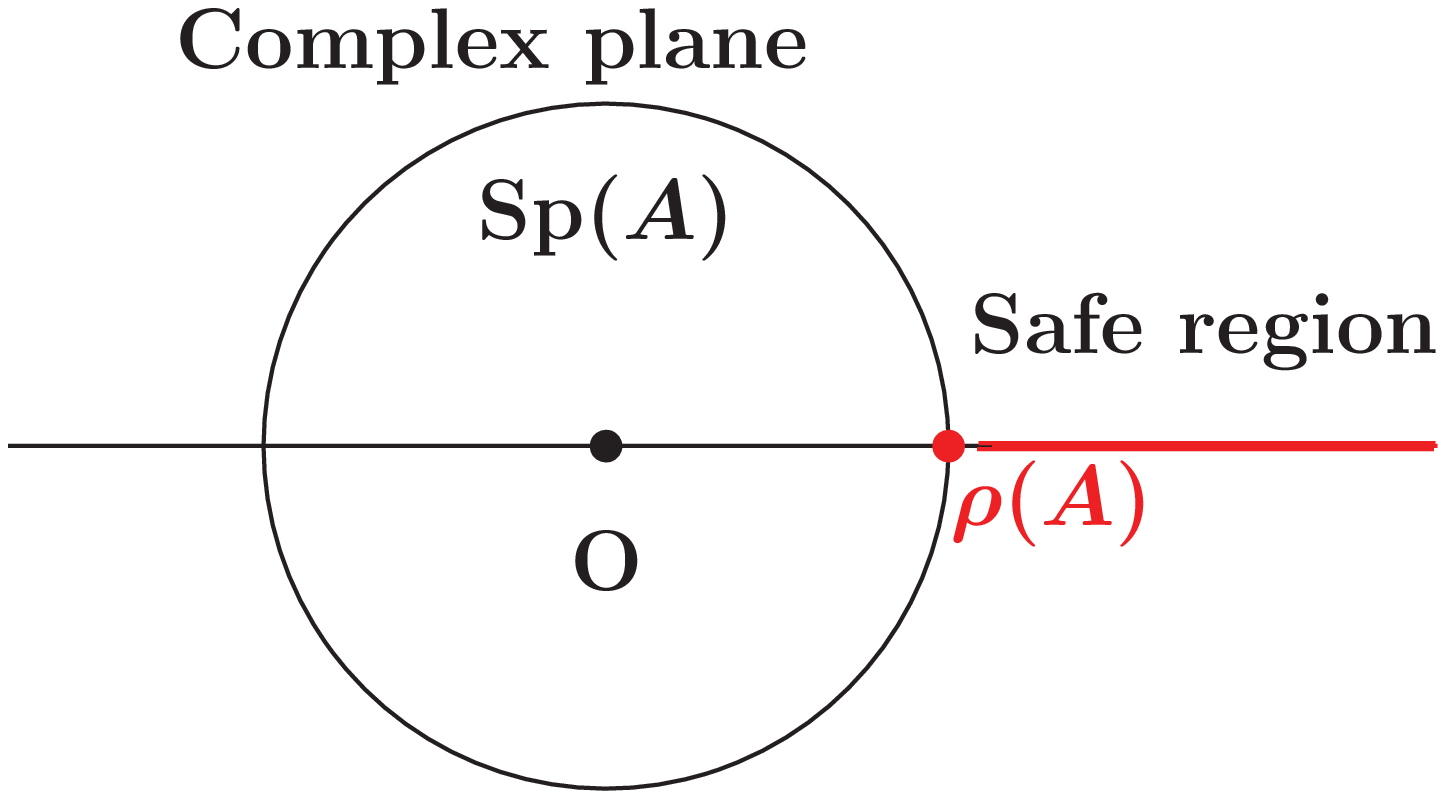}}\qd
{\includegraphics[width=6.25cm,height=3.5cm]{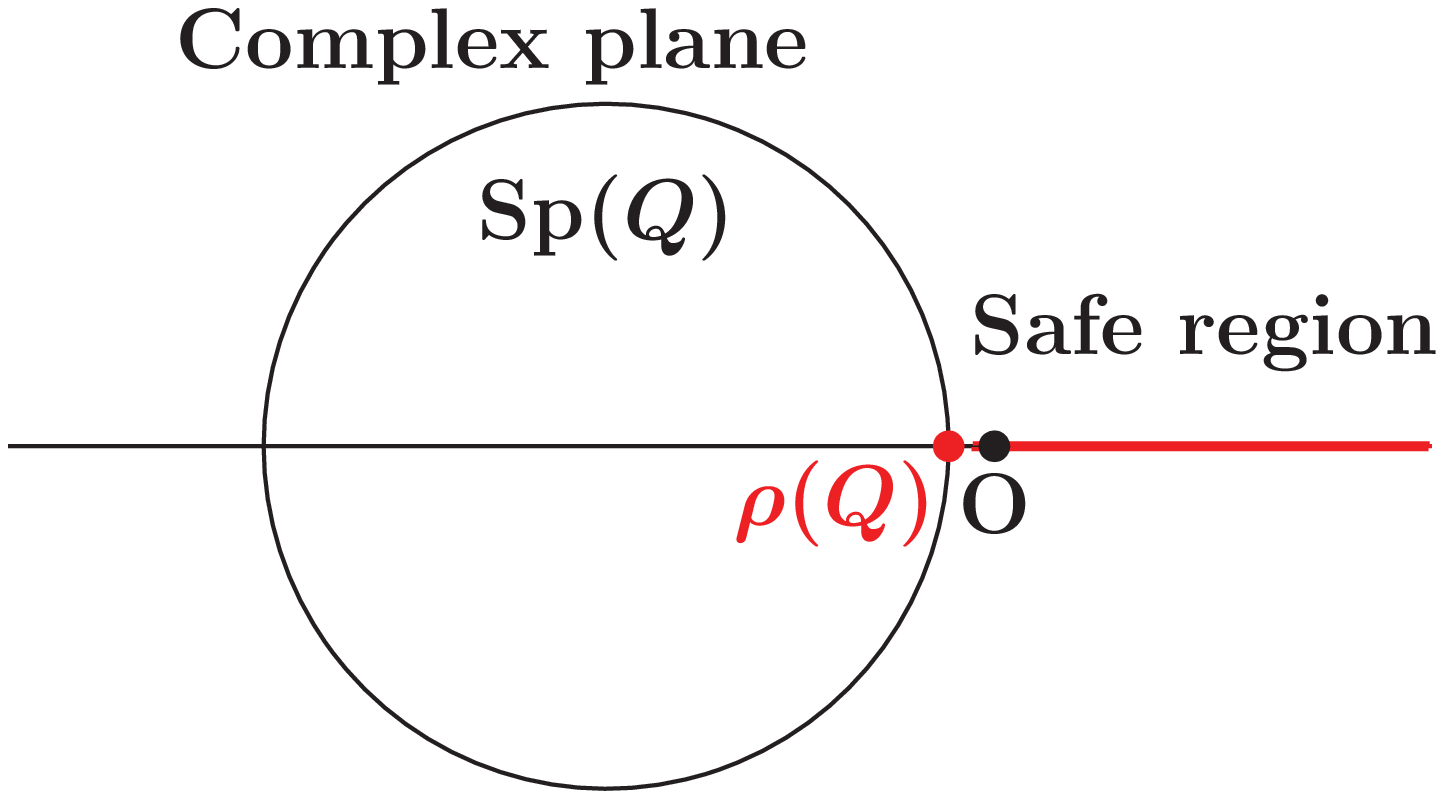}}\\
Figure 3: Safe region in complex plane.
\end{center}}

\nnd For the matrix $Q$ used in this paper, since $\rho(Q)<0$, its spectrum
$Sp(Q)$ is located on the left-hand side of the origin. Then, one
can simply choose $z_0=0$ as an initial. See Figure 3.

Having these idea in mind, we can now state two of our global algorithms.
Each of them uses the same initials:
$$v_0=\text{\rm uniform distribution},\qqd z_0=\max_{0\le i\le N}\frac{Av_0}{v_0}(i),$$
where for two vectors $f$ and $g$, $(f/g)(i)=f_i/g_i$.
\medskip

\nnd {\bf Algorithm 1}\;{\rm(Specific Rayleigh quotient iteration)}\qd
{\cms At step $k\ge 1$, for given $v:=v_{k-1}$ and $z:=z_{k-1}$, let
$w$ solve the equation
$$(z I- A) w=v.$$
Set
$v_k={w}/{\|w\|}$ and let $z_k=v_k^* A v_k.$}

This algorithm goes back to \rf{cmf16}{\S 4.1 with Choice I}.
\medskip

\nnd {\bf Algorithm 2}\;{\rm(Shifted inverse iteration)}\qd{\cms Everything is the same as in Algorithm 1, except redefine $z_k$ as follows:
$$z_k=\max_{0\le i\le N}\frac{A v_k}{v_k}(i)$$
for $k\ge 1$ (or equivalently, $k\ge 0$).}

The comparison of these algorithms is the following: with unknown small probability, Algorithm 1 is less safe than Algorithm 2, but the former one
has a faster convergence speed than the latter one with possibility 1/5 for instance.
A refined combination of the above two algorithms is presented in
\rf{cmf17c}{\S 2}, say Algorithm $4_2$ for instance.

With the worrying on the safety and convergence speed in mind, we examine
two examples which are non-symmetric.

The first example below is a lower triangular plus the upper-diagonal.
It is far away from the tridiagonal one, we want to see what can be happened.
\xmp{\rm(\rf{cmf17c}{Example 7})}\qd{\cms Let
\be Q\!=\!\begin{pmatrix}
-1 & 1 &0& 0&\cdots\cdots & 0 &0\\
a_1 & -a_1\!-\!2 &2&0& \cdots\cdots & 0 &0\\
a_2 &0\;\; & -a_2\!-\!3&3& \cdots\cdots & 0 &0\\
\vdots& \vdots &\vdots&\vdots &\cdots\cdots &N-1 &0\\
a_{N-1} & 0 &0& 0&\cdots &\!\!\!\! -a_{N-1}\!-\!N & N\\
a_N & 0 &0& 0&\cdots\cdots & 0 &\!\!-a_N\!-\!N\!-\!1
\end{pmatrix}.  \lb{04-1}\de
For this matrix, we have computed several cases:
$$a_k=1/(k+1),\;\; a_k\equiv 1,\;\; a_k=k,\;\; a_k=k^2.$$
Among them, the first one is the hardest and is hence presented below.

For different $N$,
the outputs of our algorithm are given in Table 4.   
\begin{center}{\bf Table 4}. The outputs for different $N$ by our algorithm \end{center}
\vspace{-0.4truecm}
$$\begin{tabular}{c|c|c|c|c|c|c}
\hline
\pmb{$N\!+\!1$}  & \pmb{${\displaystyle {z_1}}$} & \pmb{$z_2$} &\pmb{$z_3$} &\pmb{$z_4$} &\pmb{$z_5$} &\pmb{$z_6$}\\
 \hline\hline
$8$ & $ 0.276727$  & $0.427307$ & $0.451902$ &0.452339 & ${  } $ & ${  } $\\
 \hline
$16$ & $ 0.222132$  & $ 0.367827$ & $0.399959$ & 0.400910& ${  } $ &{  }\\
 \hline
$32$ & $ 0.187826$  & $0.329646$ & $0.370364$ & $0.372308$ & 0.372311 &{  }\\
 \hline
$50$ & $ 0.171657$  & $0.311197$ & $0.357814$ & ${0.360776} $ & 0.360784 &{  }\\
 \hline
$100$ & $0.152106$  & $0.287996$ & $0.343847$ & $0.349166$ & 0.349197 &{  }\\
 \hline
$500$ & $0.121403$  & $0.247450$ & $0.321751$ & $0.336811$ &0.337186 &{ } \\
 \hline
$1000$ & $0.111879$  & $0.233257$ & $0.313274$ & $0.334155 $ &0.335009 & 0.335010 \\
 \hline
$5000$ & $ 0.0947429$  & $0.205212$ & $0.293025$ & $0.328961$ & 0.332609 & 0.332635\\
 \hline
$10^4$ & $0.0888963$  & $0.194859$ & $0.284064$ & $0.326285$ & 0.332113 & 0.332188\\
 \hline
\end{tabular}$$
}
\dexmp

The next example is upper triangular plus lower-diagonal. It
is motivated from the classical branching process.
Denote by $(p_k: k\ge 0)$ a given probability measure with $p_1=0$. Let
$$Q\!=\!\left(\!\begin{array}{ccccccc}
-1\!\!\! & p_2 & p_3 &p_4&\cdots\cdots &p_{N-1} &\sum_{k\ge N}p_k \\
2 p_0 &\!\!\! -2 &2 p_2&2p_3 &\cdots\cdots &2 p_{N-2} & 2\sum_{k\ge N-1}p_k \\
0& 3p_0 &\!\!\! -3 & 3p_2&\cdots &3 p_{N-3} & 3\sum_{k\ge N-2}p_k\\
\vdots &\vdots &\vdots &\ddots&\ddots &\ddots &\ddots \\
\vdots &\vdots &\vdots &\ddots&\ddots &-(N\!-\!1) & (N\!-\!1)\sum_{k\ge 2}p_k\\
0& 0& 0& 0 &\cdots \cdots&\; Np_0\;\; & - Np_0
\end{array}\!\!\right)\!.$$
The matrix is defined on $E:=\{1, 2, \ldots, N\}$. Set
$M_1=\sum_{k\in E}k p_k$. When $N=\infty$, it is subcritical iff
$M_1<1$, to which the maximal eigenvalue should be positive. Otherwise, the
convergence rate should be zero.

Now, we fix
$$p_0=\az/2,\; p_1=0,\; p_2=(2-\az)/2^2,\; \ldots p_n=(2-\az)/2^n,\cdots,
\qqd \az\in (0, 2).$$
Then $M_1=3(2-\az)/2$ and hence we are in the subcritical case iff $\az\in (4/3, 2)$.

\xmp{\rm(\rf{cmf17c}{Example 9})}\qd{\cms Set $\az=7/4$. We want to know how fast the
local ($N<\infty$) maximal eigenvalue becomes stable (i.e., close enough to the
converge rate at $N=\infty$).
Up to $N=10^{4}$, the steps of the iterations we need are no more
than $6$. To quicken the convergence, we adopt an improved algorithm.
Then the outputs of the approximation of
the minimal eigenvalue of $-Q$ for different $N$ are given in Table 5.
\begin{center}{\bf Table 5}. The outputs in the subcritical case\end{center}
\vspace{-0.4truecm}
$$\begin{tabular}{c|c|c|c|c}
\hline
$N$  & \pmb{${\displaystyle {z_1}}$} & \pmb{$z_2$} &\pmb{$z_3$} &\pmb{$z_4$} \\
 \hline\hline
$8$ & $ 0.637800$  & $0.638153$ & ${  } $  & ${  } $\\
 \hline
$16$ & $ 0.621430$  & $ 0.625490$ & $0.625539$ & ${  } $\\
 \hline
$50$ & $ 0.609976$  & $ 0.624052$ & $0.624997$ & ${0.625000} $\\
 \hline
$100$ & $0.606948$  & $0.623377$ & $0.624991$ & $0.625000 $\\
 \hline
$500$ & $ 0.604409$  & $ 0.622116$ & $0.624962$ & ${0.625000} $\\
 \hline
$1000$ & $0.604082$  & $0.621688$ & $0.624944$ & $0.625000 $\\
 \hline
$5000$ & $0.603817$  & $0.620838$ & $0.62489$ & $0.625000$\\
 \hline
$10^4$ & $0.603784$  & $0.620511$ & $0.624861$ & $0.625000$\\
 \hline
\end{tabular}$$
The computation in each case costs no more than one minute.
Besides, starting from $N=50$, the final outputs are all
the same: $0.625$, which then can be regarded as a very good
approximation of $\lz_{\min}(-Q)$ at infinity $N=\infty$.
}\dexmp

It is the position to compare our global algorithm with that given
in the last section. At the first look, here in the two examples above,
we need about 6 iterations, double of the ones given in the last section.
Note that for the initials of the algorithm in the last section, we
need solve three additional linear equations, which are more or less the
same as three additional iterations. Hence the efficiency of these two
algorithms are very close to each other. Actually, the computation time
used for the algorithm in the last section is much more than the new one
here.

It is quite surprising that our new algorithms work for a much general
class of matrices, out of the scope of \ct{cmf16}. Here we consider the
maximal eigenpair only.

The example below allows partially negative off-diagonal elements.
\xmp{\rm(\rf{nout08}{Example (7)}, \rf{cmf17c}{Example 12})} \qd{\cms Let
$$A =
\begin{pmatrix} -1& 8& -1\\
8& 8& 8\\
-1& 8& 8
\end{pmatrix}.$$
Then
The eigenvalues of $A$ are as follows.
$$17.5124, \qd -7.4675,\qd 4.95513.$$
The corresponding maximal eigenvector is
$$(0.486078,\; 1.24981,\; 1)^*$$
which is positive.

Started at $z_0=24$, the outputs of our algorithms are given in Table 6.
\begin{center}{\bf Table 6}. The outputs for a matrix with more negative elements\end{center}
\vspace{-0.3truecm}
$$\begin{tabular}{c|c|c}
\hline
\pmb{$n$}  & \pmb{$z_n$: Algorithm 1} & \pmb{$z_n$: Algorithm 2} \\
 \hline\hline
$1$ & $ 17.3772$  & $18.5316$ \\
 \hline
$2$ & $ 17.5124$  & $ 17.5416$ \\
 \hline
 $3$ & ${   }$  & $ 17.5124$ \\
 \hline
\end{tabular}$$
}
\dexmp

Furthermore, we can even consider some complex matrices.

\xmp{\rm(\rf{nout12}{Example 2.1}, \rf{cmf17c}{Example 15})}\qd{\cms Let
$$A=
\begin{pmatrix}
0.75 - 1.125\, i \; &\; 0.5882 - 0.1471\, i \; &\;
  1.0735 + 1.4191\, i\\
   -0.5 -  i\; &\; 2.1765 + 0.7059\, i\; &\;
  2.1471 - 0.4118\, i\\
   2.75 - 0.125\, i \;&\;
  0.5882 - 0.1471\, i\; &\; -0.9265 + 0.4191\, i
\end{pmatrix},$$
where the coefficients are all accurate, to four decimal digits.
Then $A$ has eigenvalues
$$3,\qd -2 - i,\qd 1 + i$$
with maximal eigenvector
$$(0.408237,\;\; 0.816507,\;\; 0.408237)^*.$$
The outputs $(y_n)$ (but not $(z_n)$) of \rf{cmf17c}{Algorithm 14}, a variant of Algorithm 2, are as follows.
\begin{center}{\bf Table 7}. The outputs for a complex matrix\end{center}
\vspace{-0.6truecm}
$$\begin{tabular}{c|c|c}
\hline
\pmb{$y_1$}  & \pmb{$y_2$} & \pmb{$y_3$} \\
 \hline\hline
$3.03949 - 0.0451599\, i$ & $ 3.00471 - 0.0015769\,i$  & $3$ \\
 \hline
\end{tabular}$$
}\dexmp

We mention that a simple sufficient condition for the use of our algorithms is the following:
\be \text{\rm Re}(A^n)>0\text{ for large enough }n, \text{\rm\; up to a shift }m I.\de
Then we have the Perron--Frobenius property: there exists the maximal eigenvalue
$\rho(A)>0$ having simple left- and right-eigenvectors.

Hopefully, the reader would now be accept the use of ``global'' here for our new algorithms. They are very much efficient indeed. One may ask about the
convergence speed of the algorithms. Even though we do not have a
universal estimate for each model in such a general setup, it is
known however that the shifted inverse algorithm is a fast cubic one, and hence should
be fast enough in practice. This explains the reason why our algorithms are fast enough in the general setup. Certainly, in the tridiagonal dominate case, one can use the algorithms presented in the previous sections. Especially, in the tridiagonal situation, we have analytically basic estimates which guarantee the efficiency of the algorithms. See \ct{cmf17a} for a long way to reach the present level.

When talking about the eigenvalues, the first reaction for many people (at least for me, 30 years ago) is that well, we have known a great deal about the subject. However, it is not the trues. One may ask himself that for eigenvalues, how large matrix have you computed by hand? As far as I know, $2\times 2$ only in analytic computation by hand. It is not so easy to compute them for a $3\times 3$ matrix, except using computer. Even I have worked on the topic for about 30 years, I have not been brave enough to compute the maximal eigenvector, we use its mimic only to estimate the maximal eigenvalue (or more generally the first nontrivial eigenvalue). The first paper I wrote on the numerical computation is \ct{cmf16}. It is known that the most
algorithms in computational mathematics are local, the Newton algorithm (which is a quadratic algorithm) for instance. Hence, our global algorithms are somehow unusual.

About three years ago, I heard a lecture that dealt with a circuit board optimization problem. The author uses the Newton method. I said it was too dangerous and could fall into the trap. The speaker answered me that yes, it is dangerous, but no one in the world can solve this problem. Can we try annealing algorithm? I asked. He replied that it was too slow. We all know that in the global optimization, a big problem (not yet cracked) is how to escape from the local traps. The story we are talking about today seems to have opened a small hole for algorithms and optimization problems, and perhaps you will be here to create a new field.

\medskip

\nnd{\bf Acknowledgments}. {\small
This paper is based on a series of talks:
Central South U (2017/6),
2017 IMS-China, ICSP at Guangxi U for Nationalities (2017/6),
Summer School on Stochastic Processes at BNU (2017/7),
the 9th Summer Camp for Excellent College Students at BNU (2017/7),
Sichun U (2017/7),
the 12th International Conference on Queueing Theory and Network Applications at Yanshan U (2017/8),
the 2nd Sino-Russian Seminar on Asymptotic Methods in Probability Theory and Mathematical Statistics \& the 10th Probability Limit Theory and Statistic Large Sample Theory Seminar at Northeast Normal U (2017/9),
Workshop on Stochastic Analysis and Statistical Physics at AMSS of CAS (2017//11),
Yunnan U (2017/11).
The author thanks professors
Zhen-Ting Hou, Zai-Ming Liu,
Zhen-Qing Chen, Elton P. Hsu, Jing Yang,
Xiao-Jing Xu,
An-Min Li, Lian-Gang Peng,
Qian-Lin Li,
Zhi-Dong Bai, Ning-Zhong Shi, Jian-Hua Guo,
Zheng-Yan Lin,
Zhi-Ming Ma and C. Newman et al,
and Nian-Sheng Tang
for their invitations and hospitality. The author also thanks Ms Jing-Yu Ma
for the help in editing the paper.
Research supported in part by
         National Natural Science Foundation of China
       (Grant Nos. 11626245, 11771046),
         the ``985'' project from the Ministry of Education in China,
and the Project Funded by the Priority Academic Program Development of
Jiangsu Higher Education Institutions.
}

\vspace{-0.25truecm}

\nnd {\small
Mu-Fa Chen\\
School of Mathematical Sciences, Beijing Normal University,
Laboratory of Mathematics and Complex Systems (Beijing Normal University),
Ministry of Education, Beijing 100875,
    The People's Republic of China.\newline E-mail: mfchen@bnu.edu.cn\newline Home page:
    http://math0.bnu.edu.cn/\~{}chenmf/main$\_$eng.htm
}
\end{document}